\documentstyle{article}     

\input amssym.def
\input amssym.tex

\newtheorem{theorem}{Theorem}

\newtheorem{proposition}{Proposition}
\newtheorem{lemma}{Lemma}

\newenvironment{definition}
{\smallskip\noindent{\bf Definition\/}:}{\smallskip\par}

\newenvironment{example}
{\smallskip\noindent{\bf Example\/}.}{\smallskip\par}
\newenvironment{remark}
{\smallskip\noindent{\bf Remark\/}.}{\smallskip\par}
\newenvironment{remarks}
{\smallskip\noindent{\bf Remarks\/}.}{\smallskip\par}
\newenvironment{proof}
{\noindent{\bf Proof\/}.}{{ $\Box$}\smallskip\par}
\newenvironment{Proof}
{\noindent{\bf Proof\/}}{{ $\Box$}\smallskip\par}

\newcommand{\CC}{{\Bbb C}}

\newcommand{\RR}{{\Bbb R}}
\newcommand{\ZZ}{{\Bbb Z}}

\newcommand{\calO}{{\cal O}}

\newcommand{\calC}{{\cal C}}

\newcommand{\eps}{\varepsilon}

\title{Indices of 1-forms on an isolated complete intersection singularity}
\author{W.Ebeling and S.M.Gusein-Zade
\thanks{Partially supported by the DFG-programme ''Global methods in
complex geometry'', grants RFBR-- , INTAS--00-259.}
}
\date{}

\begin{document}

\maketitle


\section*{Introduction (\rm 1-forms versus vector fields)}

There is a number of papers devoted to the definition and the
calculation of the index of an analytic vector field on a real analytic
variety with an isolated singular point (\cite{ASV}, \cite{EGZ},
\cite{GMM1}, \cite{GMM2}). In \cite{GSV}
there was defined a notion of the index of a holomorphic vector
field on an isolated complete intersection singularity ({\bf icis}).
One can say that the reason for inventing this notion is the hope
that it can be used for calculation of the index of a vector field
on a real {\bf icis}.

We offer a different approach. Instead of considering vector
fields on a variety we consider 1-forms. One can define the
notions of the index of a real 1-form on a germ of a real
analytic variety with an isolated singular point and of a
holomorphic 1-form on a (complex) {\bf icis}. (For short
sometimes we shall refer to these two settings as to the
real and the complex indices respectively.) To a vector field
on a real variety $(V, 0)\subset (\RR^n, 0)$ (or on a complex
analytic variety $(V, 0)\subset (\CC^n, 0)$) one can associate
a 1-form on it (dependent on the choice of coordinates in
$(\RR^n, 0)$ or in $(\CC^n, 0)$). If the vector field has an
algebraically isolated singular point on $(V,0)$ (in the complex
setting this simply means isolatedness), then, for a generic
choice of coordinates, the corresponding 1-form has an algebraically
isolated singular point as well. This correspondence does not
work in the other direction. Moreover the index of a vector field on
a real analytic variety with an isolated singular point coincides
with that of the corresponding 1-form. The notion of the index
of a holomorphic 1-form on an {\bf icis} is somewhat more natural than
that of a holomorphic vector field: in some sense it is ``more
complex analytic'' (does not use the complex conjugation for the
definition) and ``more geometric'' (uses only objects of the
same tensor type). Moreover the index of an isolated singular
point of a holomorphic 1-form on a complex {\bf icis} can be described as
the dimension of an appropriate algebra. Finally, the real
index of a 1-form on a real {\bf icis} (with an algebraically
isolated singular point) plus-minus the Euler characteristic
of a (real) smoothing of the {\bf icis} can be expressed in terms of the
signature of a (nondegenerate) quadratic form on a space,
the dimension of which is equal to the (complex) index of
the corresponding complexification.

The idea to consider indices of 1-forms
instead of indices of vector fields (in some situations) was first
formulated by V.I.Arnold (\cite{Arnold79}). We are thankful to E.Looijenga,
M.Merle and J.Brian\c{c}on for useful discussions. In particular the comments
of E.Looijenga helped us to prove Theorem~\ref{theoremDim}.

\section{Indices of real 1-forms on singular varieties}

A {\em manifold with isolated singularities} is a topological
space $M$ which has the structure of a smooth (say, $C^\infty$-) manifold
outside of a discrete set $S$ (the {\em set of singular points} of $M$).
A {\em diffeomorphism} between two such manifolds is a homeomorphism
which sends the set of singular points onto the set of singular points and
is a diffeomorphism outside of them. We say that $M$ has a {\em
cone-like singularity} at a (singular) point $P\in S$ if there exists a
neighbourhood of the point $P$ diffeomorphic to the cone
$CW_P=(I\times W_P)/(\{0\}\times W_P)$ ($I=[0, 1]$) over a smooth
manifold $W_P$ ($W_P$ is called the {\em link} of the point $P$). In what
follows we assume all manifolds to have only cone-like singularities.
In \cite{EGZ} we discussed the notion of a vector field on a manifold with
isolated singularities and the notion of the index of its singular point.
Here we adapt the corresponding definitions to the case of 1-forms.
A (smooth or continuous) {\em 1-form} on a manifold $M$ with isolated
singularities is a (smooth or continuous) 1-form on the set
$M\setminus S$ of regular points of $M$. The {\em set of singular points}
$S_\omega$ of a 1-form $\omega$ on a (singular) manifold $M$ is the union
of the set of usual singular points of $\omega$ on $M\setminus S$ (i.e.,
points at which $\omega$ tends to zero) and of the set $S$ of singular
points of $M$ itself.

For an isolated {\em usual} singular point $P$ of a 1-form $\omega$
there is defined its index $\mbox{ind}_P\, \omega$ (the degree of the map
$\omega/\Vert \omega\Vert:\partial B\to S^{n-1}$ of the boundary of a
small ball $B$ centred at the point $P$ in a coordinate neighbourhood
of $P$ into the unit sphere in the dual space; $n=\mbox{dim}\,M$).
If the manifold $M$ is closed (i.e., compact, without boundary),
has no singularities ($S=\emptyset$),
and the 1-form $\omega$ on $M$ has only isolated singularities, then
\begin{equation}\label{eq1}
\sum_{P\in S_\omega}\mbox{ind}_P\, \omega=\chi(M)
\end{equation}
($\chi(M)$ is the Euler characteristic of $M$).

Let $(M, P)$ be a cone-like singularity (i.e., a germ of a manifold with
such a singular point) and let $\omega$ be a 1-form defined on an open
neighbourhood $U$ of the point $P$. Suppose that $\omega$ has no singular
points on $U\setminus\{P\}$. Let $V$ be a closed cone--like neighbourhood of
$P$ in $U$ ($V\cong CW_P$, $V\subset U$). On the cone $CW_P=(I\times W_P)
/(\{0\}\times W_P)$ ($I=[0, 1]$) there is defined a natural 1-form $dt$
($t$ is the coordinate on $I$). Let $\omega_{rad}$ be the corresponding
1-form on $V$. Let $\widetilde\omega$ be a smooth 1-form on $U$ which
coincides with $\omega$ near the boundary $\partial U$ of
the neighbourhood $U$ and with $\omega_{rad}$ on $V$ and has only isolated
singular points.

\begin{definition} The {\em index} $\mbox{ind}_P\,\omega$ of the 1-form
$\omega$ at the point $P$ is equal to
$$
1+\sum_{Q\in S_{\widetilde\omega}\setminus\{P\}}\mbox{ind}_Q\,\widetilde\omega
$$
(the sum is over all singular points $Q$ of the 1-form
$\widetilde\omega$ except $P$ itself).
\end{definition}

For a cone-like singularity at a point $P\in S$, the link $W_P$ and thus
the cone structure of a neighbourhood are, generally speaking, not
well-defined (cones over different manifolds may be {\em locally}
diffeomorphic). However it is not difficult to show that the index
$\mbox{ind}_P\, \omega$ does not depend on the choice of the cone structure
on a neighbourhood and on the choice of the 1-form $\widetilde\omega$.

\begin{proposition}\label{prop1}
For a 1-form $\omega$ with isolated singular points on a closed manifold $M$
with isolated singularities, the relation {\rm(\ref{eq1})} holds.
\end{proposition}

\begin{definition} One says that a singular point $P$ of a manifold $M$
(locally diffeomorphic to the cone $CW_P$ over a manifold $W_P$) is
{\em smoothable} if $W_P$ is the boundary of a smooth compact manifold
$\widetilde V_P$.
\end{definition}

In what follows we shall call $\widetilde V_P$ a smoothing of $(V, P)$.
The class of smoothable singularities includes, in particular, the class of
(real) isolated complete intersection singularities. For such a singularity,
there is a distinguished cone-like structure on its neighbourhood.

Let $(M, P)$ be a smoothable singularity (i.e., a germ of a manifold
with such a singular point) and let $\omega$ be a 1-form on $(M, P)$
with an isolated singular point at $P$. Let $V=CW_P$ be a closed
cone-like neighbourhood of the point $P$; $\omega$ is supposed to
have no singular points on $V\setminus\{P\}$. Let the link $W_P$ of the
point $P$ be the boundary of a compact manifold $\widetilde V_P$.
Identifying $\partial\widetilde V_P=W_P$ with $W_P\times\{1/2\}$
and using a smoothing one can consider the union
$\widetilde V_P\cup_{W_P}(W_P\times[1/2, 1])$ of $\widetilde V_P$ and
$W_P\times[1/2, 1]\subset CW_P$ glued along the common boundary as
a smooth manifold (with the boundary $W_P\times\{1\}$). The restriction
of the 1-form $\omega$ to $W_P\times[1/2, 1]\subset CW_P=V$ can
be extended to a smooth 1-form $\widetilde\omega$ on
$\widetilde V_P\cup_{W_P}(W_P\times[1/2, 1])$ with isolated singular points.

\begin{proposition}\label{prop2}
The index ${\rm ind}_P\, \omega$ of the 1-form $\omega$ at the point $P$
is equal to
$$
\sum_{Q\in S_{\widetilde \omega}}{\rm ind}_Q\, \widetilde\omega -
\chi(\widetilde V_P)+1
$$
(the sum is over all singular points of $\widetilde\omega$ on
$\widetilde V_P$).
\end{proposition}

Let $(V, 0)\subset(\RR^n, 0)$ be a real $(n-k)$-dimensional (from the
topological point of view as well) variety with
an isolated singularity at the origin.
Let $X$ be an analytic vector field on $(V, 0)$, that is the restriction
of an analytic vector field $\sum X_i\frac{\partial\ \ }{\partial x_i}$
(which we shall denote by $X$ as well) defined on a neighbourhood of
the origin in $\RR^n$ and tangent to the variety $V$ (outside of the
origin). Suppose that the origin is an isolated singular point
of the vector field $X$ on $V$, i.e., $X$ has no zeros on $V$ outside of
the origin (in a neighbourhood of it).
In this situation the index $\mbox{ind}_0\, X$ of the
vector field $X$ at the origin is defined (see \cite{EGZ}). Let $\omega$
be the 1-form $\sum X_i dx_i$. The 1-form $\omega$ on $V$ has an isolated
singular point at the origin as well. Moreover $\mbox{ind}_0 \,\omega=
\mbox{ind}_0 \, X$. Thus in this case the problem of calculating the index
of a vector field can be reduced to the problem of calculating the index
of a 1-form. This correspondence does not work in the other direction:
for a 1-form $\omega=\sum A_i dx_i$, the vector field $\sum A_i
\frac{\partial\ \ }{\partial x_i}$, generally speaking, is not tangent
to the variety $V$.

\begin{remark}
The described 1-form $\omega$ corresponding to a vector field $X$
on $V$ depends on the choice of the coordinates on $(\RR^n, 0)$.
\end{remark}

One can hope to get an algebraic formula for the index
of a vector field or of a 1-form on a singular variety
$(V, 0)\subset (\RR^n,0)$ with an isolated singular point at the origin
in the spirit of the Eisenbud--Levine--Khimshiashvili one
(\cite{EL77}, \cite{Kh}) only if
the origin is an algebraically isolated singular point of the vector
field or of the 1-form. This means that the complexification of the
vector field (or of the 1-form) on the complexification
$(V_\CC,0)\subset(\CC^n, 0)$ of the
variety $(V, 0)$ has no zeros on $V_\CC$ outside of the origin (in a
neighbourhood of it; one supposes that $V$ has an algebraically isolated
singular point at the origin itself, i.e., $V_\CC$ has an isolated singular
point at the origin). Suppose that the vector field $X=\sum X_i
\frac{\partial\ \ }{\partial x_i}$ has an algebraically isolated
singular point at the origin. In this case the 1-form $\omega =
\sum X_i dx_i$ on $V_\CC$ may have a nonisolated singular point at
it. For instance, this takes place for $V=\{(x, y, z)\in \RR^3:
x^2+y^2-z^2 = 0\}$, $X=x\frac{\partial\ }{\partial x} +
y\frac{\partial\ }{\partial y} + z\frac{\partial\ }{\partial z}$.
In such case one hardly can hope to have an algebraic formula for
the index $\mbox{ind}_0\, \omega$ of the 1-form $\omega$.
The following statement helps to avoid this difficulty.

\begin{lemma}
Let  $X=\sum X_i\frac{\partial\ \ }{\partial x_i}$ be a vector
field tangent to $V$ which has an algebraically isolated (on $V$) singular
point at the origin. Then after a possible analytic change of coordinates
in $(\RR^n, 0)$ (in fact after a generic one), the associated 1-form
$\omega=\sum X_i dx_i$ on $V$ has an algebraically isolated singular
point at the origin as well.
\end{lemma}

\begin{proof}
Consider the subset $\Xi$ of the space $J^1(V_\CC\setminus\{0\}, \CC^n)$
of 1-jets of maps from $V_\CC\setminus\{0\}$ to $\CC^n$ which consists of
jets $(F(x), dF(x))$ ($x\in V_\CC\setminus\{0\}$, $F(x)\in\CC^n$, $dF(x):
T_xV_\CC\to T_{F(x)}\CC^n=\CC^n$) such that $dF(x)(X(x))$ is orthogonal
to ${\mbox{Im\,}}dF(x)$ in the sense of the quadratic form
$\sum\limits_{i=1}^n z_i^2$ on $\CC^n$. $\Xi$ is a submanifold
of $J^1(V\setminus\{0\}, \CC^n)$ of codimension $n-k={\mbox{dim\,}V_\CC}$.
For an immersion $F:V_\CC\setminus\{0\}\to\CC^n$ the intersection points
of the image of the jet extension $j^1F$ of the map $F$ with $\Xi$
are just those points where the 1-form corresponding to the vector
field $F_*X$ vanishes on $T_{F(x)}F(V_\CC\setminus\{0\})$. The Strong
Transversality Theorem implies that, for a generic change of coordinates,
the image of $V_\CC\setminus\{0\}$ intersects $\Sigma$ (transversally)
at isolated points. Now the fact that after the change of coordinates
the 1-form corresponding to the vector field under consideration
has an algebraically isolated singular point at the origin follows from
the Curve Selection Lemma.
\end{proof}

\section{Index of a holomorphic 1-form on an icis}

Let $f=(f_1, \ldots, f_k): (\CC^n, 0)\to (\CC^k, 0)$ be an analytic
map which defines an $(n-k)$-dimensional {\bf icis}
$V=f^{-1}(0)\subset (\CC^n, 0)$ ($f_i:(\CC^n, 0)\to (\CC, 0)$).

For a germ of a holomorphic vector field $X$ which is tangent to $V$
and has an isolated singular point on $V$ at the origin,
an index is defined in \cite{GSV}
and \cite{SS96}. We recall the definition of this index. Let $U$
be a neighbourhood of the origin in $\CC^n$ where all the functions
$f_i$ ($i=1, \ldots, k$) and the vector field $X$ are defined.
Let $S_\delta \subset U$ be a sufficiently small sphere around
the origin which intersects $V$ transversally.
Let $K = V \cap S_\delta$ be the link of the {\bf icis} $(V,0)$. Define
${\mbox{grad}}\, f_i$ by
$$
{\mbox{grad}}\, f_i =
\left(\overline{\frac{\partial f_i}{\partial x_1}}, \ldots ,
\overline{\frac{\partial f_i}{\partial x_n}}\right).
$$
One has a map
$$
(X,{\mbox{grad\,}} f_1, \ldots , {\mbox{grad\,}} f_k): K\to W_{k+1}(\CC^n),
$$
where $W_{k+1}(\CC^n)$ is the Stiefel manifold of $(k+1)$-frames in
$\CC^n$. It is well-known that $W_{k+1}(\CC^n)$ is $(2(n-k)-2)$--connected and
that $H_{2(n-k)-1}(W_{k+1}(\CC^n)) \cong
\pi_{2(n-k)-1}(W_{k+1}(\CC^n)) \cong \ZZ$ (see, e.g., \cite{H}).
On the other hand, $K$ is a smooth manifold of dimension $2(n-k)-1$
with the natural orientation as the boundary of the complex manifold
$V\setminus\{0\}$. Therefore a map from $K$ to
$W_{k+1}(\CC^n)$ has a degree. The index of the vector field $X$
at the origin is defined to be the degree of the above map.

\begin{remark}
Pay attention that one uses the complex conjugation for this definition
and the components of the discussed map are of different tensor nature.
Whereas $X$ is a vector field, ${\mbox{grad\,}}f_i$ is more similar
to a covector.
\end{remark}

Now we adapt the definition of the index of a holomorphic
vector field to the case of a holomorphic 1-form on $V$.
Let $\omega=\sum A_idx_i$ ($A_i=A_i(x)$) be a germ of a
holomorphic 1-form on $(\CC^n, 0)$ which
as a 1-form on $V$ has (at most) an isolated singular point at the
origin (thus it does not vanish on the tangent space $T_PV$ to the
variety $V$ at all points $P$ from a punctured neighbourhood
of the origin in $V$). The 1-forms $\omega$, $df_1$, \dots , $df_k$
are linearly independent for all $P \in K$. Thus one has a map
$$(\omega, df_1, \ldots , df_k) : K \to W_{k+1}(\CC^n).$$

\begin{definition}
We define the {\em index} of the
1-form $\omega$ at $0$, ${\rm ind}_{\CC,0}\, \omega$, to be the degree
of the map
$$(\omega, df_1, \ldots , df_k) : K \to W_{k+1}(\CC^n)$$
(here $W_{k+1}(\CC^n)$ is the manifold of $(k+1)$-frames in
the dual $\CC^n$).
\end{definition}

\begin{remark}
The Stiefel manifold $W_{k+1}(\CC^n)$ of $(k+1)$-frames in $\CC^n$
is homotopy equivalent to the (Stiefel) manifold $\widetilde W_{k+1}(\CC^n)$
of orthonormal (with respect to the Hermitian scalar product $\sum
x_i\overline{y_j}$)
$(k+1)$-frames in $\CC^n$. The homotopy equivalence
is defined by the Gram-Schmidt process. Thus in the definition one can
substitute $W_{k+1}(\CC^n)$ by $\widetilde W_{k+1}(\CC^n)$ and the map
$$
(\omega, df_1, \ldots , df_k) : K \to W_{k+1}(\CC^n)
$$
by the corresponding map
$$
(\omega, df_1, \ldots , df_k)\widetilde{\ }: K \to \widetilde W_{k+1}(\CC^n).
$$
However $\widetilde W_{k+1}(\CC^n)$ is not a complex analytic manifold and thus
the map
$$
(\omega, df_1, \ldots , df_k)\widetilde{\ }: U\setminus\{0\} \to
\widetilde W_{k+1}(\CC^n)
$$
of a punctured neighbourhood $U\setminus\{0\}$
of the origin in $V$ is not complex analytic (in contrast with the map
$(\omega, df_1, \ldots , df_k) : U\setminus\{0\} \to W_{k+1}(\CC^n)$).
We prefer to give a ``more complex analytic'' definition, using sometimes
the map to $\widetilde W_{k+1}(\CC^n)$ for calculations.
\end{remark}

\begin{example}
Let $n=2$, $k=1$, $f_1(x,y)=x^2+y^3$ and consider the 1-form
$$
\omega=3y^2dx-2xdy.
$$
The form $\omega$ on $V=\{f_1=0\}$ has an isolated zero at the origin.
One can easily compute that the degree of the map
$$
(\omega, df_1)\widetilde{\ } : K \to \widetilde W_2(\CC^2) \cong U(2)
$$
is 6. Therefore ${\rm ind}_{\CC,0}\, \omega = 6$ (see also
Proposition~\ref{propZeros}). This 1-form corresponds (in
the described sense) to the vector field
$$
X = 3y^2 \frac{\partial}{\partial x} - 2 x\frac{\partial}{\partial y}.
$$
This vector field (by chance) is tangent to the hypersurface
$V=\{f_1=0\}$ and has an isolated singular point at the origin
on it, but its index is $0$, since the map
$$
(X, {\mbox{grad\,}}f_1)\widetilde{\ } : K \to \widetilde W_2(\CC^2) \cong U(2)
$$
maps $K$ to $SU(2)$ which is simply connected.
\end{example}

Let $B_\delta$ be the ball of radius $\delta$ around the origin in $\CC^n$
with the boundary $S_\delta$. Suppose that the functions $f_1$, \dots , $f_k$
and the 1-form $\omega$
are defined on a neighbourhood of $B_\delta$ and let $\eps=(\eps_1, \ldots ,
\eps_k) \in \CC^k$ be small enough (so that the level set
$V_\eps := f^{-1}(\eps) \cap B_\delta$ is transversal to the sphere
$S_\delta$) and such that this level set $V_\eps$ is nonsingular.

\begin{definition}
Let $M$ be a complex $m$-dimensional manifold and let $\eta$ be a
holomorphic 1-form on $M$. A zero $P \in M$ of the 1-form $\eta$ is called
{\em nondegenerate} if in local coordinates $y_1$, \dots , $y_m$
on $M$ around $P$ in which the 1-form $\eta$ is written as
$\eta = C_1 dy_1 + \cdots C_m dy_m$, the Hessian of $\eta$ at
the point $P$, i.e., the determinant of the matrix
$$
\left(
\frac{\partial C_i}{\partial y_j}(P)
\right)_{i,j=1,\ldots, m},
$$
is nonzero.
\end{definition}

There exists a perturbation $\widetilde\omega$ of the
1-form $\omega$ which has only nondegenerate zeros on $V_\eps$.
(In fact a generic perturbation of $\omega$ and, in particular, a
perturbation of the form $\widetilde{\omega} = \omega - \lambda \eta$
for a generic 1-form $\eta$ on $\CC^n$ (where $\lambda\ne 0$ is small enough)
possesses this property.)

\begin{proposition}\label{propZeros}
The index ${\rm ind}_{\CC, 0}\, \omega$ of the 1-form $\omega$
on the {\bf icis} $V$ at the origin
is equal to the number of zeros of $\omega$ on $V_\eps$, counted with
multiplicities. It is also equal to the number of zeros of $\widetilde{\omega}$
on $V_\eps$ for a small perturbation $\widetilde{\omega}$ of the 1-form
$\omega$ with only nondegenerate zeros on $V_\eps$.
\end{proposition}

\begin{proof}
Let $P_1$, \dots, $P_\nu$ be the zeros of the form $\omega$ on $V_\eps$.
In local coordinates $y_1$, \dots, $y_{k}$, $y_{k+1}$, \dots, $y_n$
centred at the point $P_i$ such that $V_\eps=\{y_1=\ldots=y_k=0\}$
(one can take $y_i=f_i-\varepsilon_i$ for $1\le i\le k$) let
$\omega_{\vert V_\eps}=
C_{k+1}dy_{k+1}+\ldots+C_n dy_n$. Let $B_i$ be a small open ball
centred at
the point $P_i$. The degree of the map $(\omega, df_1, \ldots , df_k):
\partial B_i \cap V_\eps\to W_{k+1}(\CC^n)$ is equal to the degree
of the map $S^{2(n-k)-1} \to S^{2(n-k)-1}$ given by
$$
x \mapsto (C_{k+1}, \ldots , C_n)/\Vert(C_{k+1}, \ldots , C_n)\Vert
$$
and thus is equal to the multiplicity $\mu_i$ of the zero $P_i$
(see, e.g., \cite{AGV}). Now consider the manifold
$$M:= V_\eps \setminus \bigcup_i B_i.$$
Since the 1-form $\omega$ has no zeros on $M$, the degree of the mapping
$$(\omega, df_1, \ldots , df_k) : \partial M \to W_{k+1}(\CC^n)$$
is equal to zero. This implies
$${\rm ind}_{\CC,0} \, \omega = \sum_{i=1}^\nu \mu_i.$$
If $\widetilde{\omega}$ is a perturbation of the 1-form $\omega$ with only
nondegenerate zeros on $V_\eps$ then it has just $\mu_i$ nondegenerate
zeros on $B_i\cap V_\eps$.
\end{proof}

\section{An algebraic formula for the index}

Our aim in this section is to derive the following algebraic formula for
the complex index.

\begin{theorem} \label{theoremDim}
Let $\calO_{\CC^n,0}$ be the algebra of germs of holomorphic functions
at the origin in $\CC^n$ and let
$I \subset \calO_{\CC^n,0}$ be the ideal generated by $f_1, \ldots,
f_k$ and the
$(k+1) \times (k+1)$-minors of the matrix
$$
\left( \begin{array}{ccc} \frac{\partial f_1}{\partial x_1} & \cdots &
\frac{\partial f_1}{\partial x_n} \\
\vdots & \ddots & \vdots \\
\frac{\partial f_k}{\partial x_1} & \cdots & \frac{\partial f_k}{\partial
x_n}\\
A_1 & \cdots & A_n
\end{array} \right).
$$
Then
$${\rm ind}_{\CC,0}\, \omega = \dim_\CC \calO_{\CC^n,0}/I.$$
\end{theorem}

\begin{remark}
A corresponding result for the case when the 1-form $\omega$ is the
differential $df_{k+1}$ of a function $f_{k+1}$ was first proven by
G.-M.Greuel in \cite[Lemma~1.9]{G}. G.-M.Greuel has informed us that his
arguments in that paper can be adapted to prove Theorem~\ref{theoremDim}
as well.
\end{remark}

The proof will rely on a basic fact which we now want to formulate.

Let $F:(\CC^n\times\CC^M, 0)\to (\CC^k\times\CC^M, 0)$ be a versal
deformation of $f: (\CC^n, 0) \to (\CC^k,0)$, let $N:=n+M$, $K:=k+M$.
By the same symbol, $F$, we denote a representative $U \to \CC^K$
of this deformation defined in a small open neighbourhood $U \subset \CC^N$
of the origin. Let
$A_i=0$ for $i=n+1, \ldots , N$ and let $\eta = \sum_{i=1}^N B_i dx_i$
be a 1-form such that for a regular value $s \in F(U)$ of $F$ and
for $\lambda \in W$, $\lambda \neq 0$, where $W \subset \CC$ is a
suitable small open neighbourhood of the origin, the form
$\omega - \lambda \eta$ has only isolated zeros on $F^{-1}(s)$.

Consider the matrix
$$\Phi = \left( \begin{array}{ccc} \frac{\partial F_1}{\partial x_1} &
\cdots &
\frac{\partial F_1}{\partial x_N} \\
\vdots & \ddots & \vdots \\
\frac{\partial F_K}{\partial x_1} & \cdots & \frac{\partial F_K}{\partial
x_N}\\
A_1 - \lambda B_1 & \cdots & A_N - \lambda B_N
\end{array} \right).$$

Let $\CC^{N+1}$ be the vector space with coordinates $x_1$, \dots, $x_N$,
$\lambda$. Let $\calC( U \times W)$ denote the ideal in
$\calO_{\CC^{N+1}}(U \times W)$ generated by the $(K+1) \times (K+1)$
minors of the matrix $\Phi$. Let $\calC \subset \calO_{\CC^{N+1}}$
denote the corresponding ideal sheaf and let $C \subset U\times W$
be the analytic space defined by $\calC$. Then $C$ consists of those points
$x \in U \times W$ where $F$ is not a submersion or $\omega - \lambda \eta$
does not have an isolated zero in $x \in F^{-1}(F(x))$. Let
$\calO_C := \calO_{\CC^{N+1}}/\calC$ be the structure sheaf of $C$.

Let $\widetilde{F} : U\times\CC \to \CC^{K+1}$ be the mapping defined by
$(x,\lambda)
\mapsto (F(x), \lambda)$.
Let $\Sigma := \widetilde{F}(C)$ and let $\pi := \widetilde{F}|_C : C \to
\Sigma$.
Then $\Sigma$ is an analytic space endowed with the structure sheaf
$\calO_\Sigma := \calO_{\CC^{K+1}}/{\cal F}_0(\pi_\ast(\calO_C))$
(restricted to
$\pi(C)$) where ${\cal F}_0$ denotes the 0-th Fitting ideal (cf.\
\cite[4.E]{Looijenga84}). The mapping $\pi : C \to \Sigma$ is a finite mapping.

Theorem~\ref{theoremDim} will follow from the next theorem.

\begin{theorem} \label{theoremFlat}
For every $x \in C$, $\calO_{C,x}$ is a flat $\calO_{\Sigma,\pi(x)}$-module.
\end{theorem}

\begin{proof} The proof follows the same lines as \cite[(4.4) and
(4.8)]{Looijenga84}.

Let $x \in C$ and $s=\pi(x)$.
The minors of the matrix $\Phi$ vanish at a point $x \in U \times
W$ if and only if the rank of the matrix $\Phi(x)$ is not maximal. The set
of $(N+1)\times (K+1)$ matrices with complex entries with rank
$< K+1$ is an affine algebraic variety of codimension
$N-K+1$ inside the set of all $(N+1) \times (K+1)$ matrices with complex
entries. Hence
$C$ has dimension $K$. This implies that ${\rm depth}\, ( \calC_x;
\calO_{\CC^{N+1},x}) =
N-K+1$. By \cite[Corollary~(2.7)]{BR64} this implies that the homological
dimension
${\rm hd}\, \calC_x$ of $\calC_x$ is equal to $N-K+1$.

By the formula of Aus\-lan\-der--Buchsbaum (see,
e.g., \cite[p.~114]{Matsumura80}) we have
$${\rm depth}\, \calO_{C,x} = N+1 - (N-K+1) = K.$$
Since also $\dim \calO_{C,x} = K$, it follows that $\calO_{C,x}$ is a
Cohen-Macaulay ring.

Since $\widetilde{F}|_C : C \to \CC^{K+1}$ is a finite mapping, $\calO_{C,x}$
is a finite $\calO_{\CC^{K+1},s}$-module. Therefore $\calO_{C,x}$ is a
Cohen-Macaulay
module over $\calO_{\CC^{K+1},s}$ and
$$
{\rm depth}_{\calO_{\CC^{K+1},s}}\, \calO_{C,x} = \dim \calO_{C,x} = K.
$$
Since $\calO_{\CC^{K+1},s}$ is an $(K+1)$-dimensional regular ring, it
follows from the
Aus\-lan\-der--Buchsbaum formula that the homological dimension of
$\calO_{C,x}$
as an
$\calO_{\CC^{K+1},s}$-module is 1. This means that there is an exact
sequence of
$\calO_{\CC^{K+1},s}$-modules
$$0 \rightarrow \calO_{\CC^{K+1},s}^q \stackrel{\alpha}{\rightarrow}
\calO_{\CC^{K+1},s}^p \rightarrow \calO_{C,x} \rightarrow 0.$$
Moreover, $q$ must be equal to $p$ and the 0-th Fitting ideal of
$\calO_{C,x}$ viewed as
an $\calO_{\CC^{K+1},s}$-module is generated by the determinant of
$\alpha$. Hence
$\Sigma$ is a hypersurface and $\calO_{\Sigma,s}$ is a Cohen-Macaulay ring,
too.

Since $\pi : C \to \Sigma$ is a finite mapping, $\calO_{C,x}$
is also a finite $\calO_{\Sigma,s}$-module. But a finitely generated
$\calO_{\Sigma,s}$-module is flat if and only if it is free (see e.g.\
\cite[(3.G) Proposition]{Matsumura80}). Therefore it suffices to show that
$\calO_{C,x}$
is  a free $\calO_{\Sigma,s}$-module. By the Aus\-lan\-der--Buchsbaum
formula we
have
$$
{\rm hd}_{\calO_{\Sigma,s}}\, \calO_{C,x} + {\rm depth}_{\calO_{\Sigma,s}}\,
\calO_{C,x} =  {\rm depth}\, \calO_{\Sigma,s}.
$$
Since ${\rm depth}_{\calO_{\Sigma,s}}\, \calO_{C,x} = {\rm depth}\,
\calO_{\Sigma,s} =
\dim
\calO_{\Sigma,s} = K$, it follows that ${\rm hd}_{\calO_{\Sigma,s}}\,
\calO_{C,x} = 0$.
But this means that
$\calO_{C,x}$ is a free $\calO_{\Sigma,s}$-module.
\end{proof}

\addvspace{3mm}

\begin{Proof} {\bf of Theorem~\ref{theoremDim}.} Consider again the
mapping $\pi : C
\to \Sigma$ and let $s \in \Sigma$. For each $x \in C(s) := \pi^{-1}(s)$, $\CC
\otimes_{\calO_{\Sigma,s}} \calO_{C,x}$ is a finite dimensional vector
space over $\CC$.
Denote its dimension by $\nu(x)$. Define
$$\nu(s) = \sum_{x \in C(s)} \nu(x).$$
By Theorem~\ref{theoremFlat} and \cite[\S5, Theorem~1]{Douady68}, $\nu(s)$
is a locally
constant function of $s$. Now for a point $s=(s',\lambda)$ where $s'$ is a
regular value
of $F$,
$\nu(s)= {\rm ind}_{\CC,0}\,
\omega$ by Proposition~\ref{propZeros}. On the other hand $\nu(0) =
\dim_\CC
\calO_{\CC^n,0}/I$.
\end{Proof}

\section{The real index as the signature of a family of quadratic forms}

Now we want to discuss the index of a real analytic 1-form
$\omega$
on a real {\bf icis} $(V, 0)=\{f_1=\ldots=f_k=0\}\subset(\RR^n, 0)$
(with an algebraically isolated singular point at the origin). The
1-form $\omega$ on $V$ is the restriction of an analytic 1-form
$\sum\limits_{i=1}^n A_1dx_i$ defined
on $\RR^n$ in a neighbourhood of the origin (we also denote this 1-form by
$\omega$).  We consider the (analytic) functions $f_1$, \dots, $f_k$
and the 1-form $\omega$ defined on a neighbourhood
of the origin in $\CC^n$ as well. Let $\delta$ be positive and
small enough so that the functions $f_1$, \dots, $f_k$ and the
1-form $\omega$ are defined on the ball $B_\delta$ of radius $\delta$
centred at the origin in $\CC^n$ and for each positive $\delta^\prime<\delta$
the variety $V_\CC=\{f_1=\ldots=f_k=0\}\subset(\CC^n, 0)$
intersects the sphere $S_{\delta^\prime}$ of radius
$\delta^\prime$ centred at the origin transversally. For
$\varepsilon=(\varepsilon_1, \ldots,\varepsilon_k)$ small enough
($\vert\varepsilon\vert \ll \delta$),
let $V_{\CC,\varepsilon}=\{f=\varepsilon\}\cap B_\delta\subset\CC^n$,
for $\varepsilon$ real,
$V_\varepsilon=V_{\CC,\varepsilon}\cap\RR^n$.
Let $\Sigma\subset(\CC^k,0)$ be the bifurcation diagram of the map
$f$ (the set of $\varepsilon\in(\CC^k,0)$
critical for $f$), $\Sigma_\RR=\Sigma\cap(\RR^k, 0)$.

By Proposition~\ref{prop2}, for real $\varepsilon$ outside of the
bifurcation diagram (i.e., $\varepsilon\in \RR^k\setminus \Sigma_\RR$),
the (real) index $\mbox{ind}_0\, \omega$ differs from
the sum of indices of zeros of the 1-form $\omega$ on the real
smooth manifold $V_\varepsilon$ by ($\chi(V_\varepsilon)-1$).
By Proposition~\ref{propZeros} the complex index $\mbox{ind}_{\CC,0}\, \omega$
counts the number of zeros of the same 1-form on the complex manifold
$V_{\CC,\varepsilon}$. The Euler characteristic $\chi(V_\varepsilon)$
of the real level manifold $V_\varepsilon$ is different for different
$\varepsilon\in\RR^k\setminus\Sigma_\RR$ (at least for $n-k$ even).
It is constant on each component of the
complement to the bifurcation diagram $\Sigma_\RR$. Thus one cannot
expect to get the (real) index $\mbox{ind}_0\, \omega$ of the 1-form $\omega$
as the signature of a nondegenerate quadratic form on a vector space of
dimension $\mbox{ind}_{\CC,0}\, \omega$ ($=\dim_\CC \calO_{\CC^n,0}/I$).
Such a signature may be equal to $\mbox{ind}_0\,
\omega+(\chi(V_\varepsilon)-1)$
and thus has to be different for different components of the complement to
the bifurcation diagram. We want to obtain a somewhat finer 
description.

\begin{theorem}\label{theo3}
There exists a family $Q_{\varepsilon}$ of quadratic
forms on the space $\CC^L$ of dimension $L=\dim_\CC \calO_{\CC^n,0}/I$
$($i.e., a family of symmetric $L\times L$-matrices$)$
defined for $\varepsilon\in \CC^k$ from
a neighbourhood of the origin and analytically dependent on
$\varepsilon$ such that: \newline
1) for $\varepsilon$ from the complement to the bifurcation diagram $\Sigma$,
the quadratic form $Q_{\varepsilon}$ is nondegenerate; \newline
2) for real $\varepsilon$, the quadratic form $($i.e., the matrix$)$
$Q_{\varepsilon}$ is real and, for real $\varepsilon$ outside of the
bifurcation diagram $($i.e., for $\varepsilon\in \RR^k\setminus\Sigma_\RR$$)$,
its signature is equal to
\begin{equation}\label{eq2}
\sum\limits_{P\in V_{\varepsilon}}{\rm ind}_P\, \omega={\rm ind}_0\, \omega+
(\chi(V_{\varepsilon})-1).
\end{equation}
\end{theorem}

\begin{proof}
It is convenient to define the family $Q_{\varepsilon}$ for a larger
space of parameters.
Let $F:(\CC^n\times\CC^M,0)\to(\CC^k\times\CC^M,0)$ be a real
(i.e., invariant with respect to the complex conjugation) versal
deformation of the $f:(\CC^n, 0)\to (\CC^k,0)$, $F(x,\varepsilon^\prime)=
(f_{\varepsilon^\prime}(x), \varepsilon^\prime)$, $f_0=f$. Here
$\varepsilon_{k+1}$, \dots, $\varepsilon_{k+M}$ are the coordinates on
$\CC^M$, $\varepsilon^\prime=(\varepsilon_{k+1}, \ldots, \varepsilon_{k+M})
\in \CC^M$; let $\hat\varepsilon=(\varepsilon_1, \ldots, \varepsilon_k,
\varepsilon_{k+1}, \ldots, \varepsilon_{k+M})\in \CC^{k+M}=\CC^k\times\CC^M$.
Let $\CC^n_\alpha$ be the $n$-dimensional affine space
with the coordinates $\alpha_1$, \dots, $\alpha_n$, let $\varkappa=
(\hat\varepsilon, \alpha)$ ($\varkappa\in\CC^{k+M+n}_\varkappa=\CC^{k+M}\times
\CC^n_\alpha$). Let us denote by the same letter, $F$, the trivial
extension $(\CC^n_x\times\CC^M_{\varepsilon^\prime}\times\CC^n_\alpha,0)\to
(\CC^k_\varepsilon\times\CC^M_{\varepsilon^\prime}\times\CC^n_\alpha,0)
=(\CC^{k+M+n}_\varkappa, 0)$ of the chosen versal deformation:
$F(x,{\varepsilon^\prime},\alpha)=
(f_{\varepsilon^\prime}(x),{\varepsilon^\prime},\alpha)$. Let $\Omega$ be
a 1-form on $(\CC^n_x\times\CC^M_{\varepsilon^\prime}\times\CC^n_\alpha,0)$
defined by
$$
\Omega=\omega-\sum\limits_{i=1}^n\alpha_i dx_i=
\sum\limits_{i=1}^n(A_i-\alpha_i)dx_i.
$$
Again by the same symbols, $\Sigma$ and $\Sigma_\RR$, we shall denote
the complex and the real bifurcation sets of the map $F$ (the real one
being the intersection of the complex one with the real space $\RR^{k+M+n}$),
$V_{\CC,\varkappa}:=F^{-1}(\varkappa)\cap B_\delta\subset\CC^{n+M+n}$,
for $\varkappa$ real,
$V_\varkappa:=V_{\CC,\varkappa}\cap\RR^{n+M+n}$ (in fact
$V_{\CC,\varkappa}$ is a subvariety of the $n$-dimensional affine
space $\CC^n_x\times\{(\varepsilon^\prime, \alpha)\}$).

For a point $P\in\CC^n_x\times\CC^M_{\varepsilon^\prime}\times\CC^n_\alpha$
let
$$
\Delta=\Delta(P) :=
\left|
\begin{array}{ccc}
\frac{\partial f_{\varepsilon^\prime 1}}{\partial x_1} &\cdots &
\frac{\partial f_{\varepsilon^\prime 1}}{\partial x_k}\\
\vdots & \ddots & \vdots \\
\frac{\partial f_{\varepsilon^\prime k}}{\partial x_1} & \cdots &
\frac{\partial f_{\varepsilon^\prime k}}{\partial x_k}
\end{array}
\right|,
$$
where $f_{\varepsilon^\prime i}$ ($i=1, \ldots, k$) are the components
of the map $f_{\varepsilon^\prime}$: $f_{\varepsilon^\prime} =
(f_{\varepsilon^\prime 1}, \ldots, f_{\varepsilon^\prime k})$.
If $P\in V_{\CC,\varkappa}$ and $\Delta(P)\ne 0$, the variables $x_{k+1}$,
\dots, $x_n$ are local coordinates on the ($(n-k)$-dimensional)
variety $V_{\CC,\varkappa}$ in a neighbourhood of the point $P$.
For a generic $\varkappa\in(\CC^{k+M+n},0)\setminus\Sigma$, for
all zeros $P$ of the restriction of the 1-form $\Omega$ to the level
manifold $V_{\CC,\varkappa}$, one has $\Delta(P)\ne 0$. Moreover,
for a generic $\varkappa\in(\CC^{k+M+n},0)\setminus\Sigma$ which does
not possess this property (such $\varkappa$'s form a subset of codimension 1),
all zeros $P$ of the restriction of the 1-form $\Omega$ to the level
manifold $V_{\CC,\varkappa}$ are simple (i.e., nondegenerate).

For a local system of coordinates $y_1$, \dots, $y_m$ on a manifold $M$
the Hessian $h=h_\omega$ of a 1-form $\omega$ equal to
$A_1dy_1+\ldots+A_mdy_m$ in these coordinates is defined as the determinant
$$
\left|
\frac{\partial A_i}{\partial y_j}
\right|_{i,j=1, \ldots , m}.
$$
(The Hessian is a function on the manifold $M$; it depends on the choice
of the coordinates.)

For $i = k+1, \ldots , n$, let
$$
m_{i}
:= \left|
\begin{array}{cccc}
\frac{\partial f_{\varepsilon^\prime 1}}{\partial x_1} &
\cdots &
\frac{\partial f_{\varepsilon^\prime 1}}{\partial x_k} &
\frac{\partial f_{\varepsilon^\prime 1}}{\partial x_i} \\
\vdots & \ddots & \vdots & \vdots \\
\frac{\partial f_{\varepsilon^\prime k}}{\partial x_1} &
\cdots &
\frac{\partial f_{\varepsilon^\prime k}}{\partial x_k} &
\frac{\partial f_{\varepsilon^\prime k}}{\partial x_i}
\\ A_1-\alpha_1 & \cdots & A_k-\alpha_k & A_i-\alpha_i
\end{array}
\right|.
$$
For a point $P\in V_{\CC, \varkappa}$ with $\Delta(P)\ne 0$,
let $h(P)$ be the Hessian of the 1-form $\Omega_{\vert V_{\CC, \varkappa}}$
in the (local) coordinates $x_{k+1}$, \dots, $x_n$.

\begin{proposition}\label{propHessian}
The Hessian $h$ of the 1-form $\Omega_{\vert V_{\CC, \varkappa}}$
in the coordinates $x_{k+1}$, \dots, $x_n$ is given by the formula:
$$
h =
\frac{1}{\Delta^{2+(n-k)}}
\left|
\begin{array}{cccc}
\Delta & \frac{\partial \Delta}{\partial x_1} & \cdots &
\frac{\partial\Delta}{\partial x_n}\\
0 & \frac{\partial f_{\varepsilon^\prime 1}}{\partial x_1} & \cdots &
\frac{\partial f_{\varepsilon^\prime 1}}{\partial x_n} \\
\vdots & \vdots & \ddots & \vdots \\
0 & \frac{\partial f_{\varepsilon^\prime k}}{\partial x_1} & \cdots &
\frac{\partial f_{\varepsilon^\prime k}}{\partial x_n} \\
m_{k+1} & \frac{\partial m_{k+1}}{\partial x_1} & \cdots &
\frac{\partial m_{k+1}}{\partial x_n}\\
\vdots & \vdots & \ddots & \vdots \\
m_n & \frac{\partial m_n}{\partial x_1} & \cdots &
\frac{\partial m_n}{\partial x_n}
\end{array}
\right|.
$$
\end{proposition}

\begin{lemma} \label{lemmaImplicit}
In the coordinates $x_{k+1}$, \dots, $x_n$ the 1-form $\Omega_{\vert
V_{\CC, \varkappa}}$
can be expressed as
$$
\Omega_{\vert V_{\CC, \varkappa}} =  \frac{m_{k+1}}{\Delta} dx_{k+1} +
\cdots + \frac{m_n}{\Delta}dx_n.
$$
\end{lemma}

\begin{proof} One has
\begin{eqnarray*}
0 & = & df_{\varepsilon^\prime 1} =
\frac{\partial f_{\varepsilon^\prime 1}}{\partial x_1} dx_1 + \cdots +
\frac{\partial f_{\varepsilon^\prime 1}}{\partial x_n} dx_n, \\
\vdots & \vdots & \vdots \\
0 & = & df_{\varepsilon^\prime k} =
\frac{\partial f_{\varepsilon^\prime k}}{\partial x_1} dx_1 + \cdots +
\frac{\partial f_{\varepsilon^\prime k}}{\partial x_n} dx_n.
\end{eqnarray*}
By Cramer's rule 
\begin{eqnarray*}
dx_1 & = & (-1)^k \frac{1}{\Delta} \left(
\frac{\partial(f_{\varepsilon^\prime 1}, \ldots,
f_{\varepsilon^\prime k})}{\partial(x_2, \ldots ,x_k, x_{k+1})} dx_{k+1} +
\cdots +
\frac{\partial(f_{\varepsilon^\prime 1},
\ldots, f_{\varepsilon^\prime k})}{\partial(x_2, \ldots , x_k, x_n)} dx_n
\right),\\
\vdots & \vdots & \vdots \\
dx_k & = & - \frac{1}{\Delta} \left( \frac{\partial(f_{\varepsilon^\prime
1}, \ldots,
f_{\varepsilon^\prime k})}{\partial(x_1, \ldots ,x_{k-1}, x_{k+1})}
dx_{k+1} + \cdots +
\frac{\partial(f_{\varepsilon^\prime 1},
\ldots, f_{\varepsilon^\prime k})}{\partial(x_1, \ldots , x_{k-1}, x_n)}
dx_n \right).
\end{eqnarray*}
Substitution in the form $\Omega$ yields
$$\Omega_{\vert V_{\CC, \varkappa}} =
 \frac{m_{k+1}}{\Delta} dx_{k+1} + \cdots + \frac{m_n}{\Delta}dx_n.$$
\end{proof}

By some abuse of notation we shall denote partial derivatives
in the coordinates $x_{k+1}$, \dots, $x_n$ of functions on the manifold
$V_{\CC,\varkappa}$ by $\frac{d\ \ }{dx_j}$ ($j=k+1,\ldots, n$)
(in order to distinguish them from partial derivatives of the
corresponding functions on $\CC^n$).
Therefore the Hessian $h$ is the determinant of the matrix
$$\left( \frac{d }{d x_j} \left(\frac{m_i}{\Delta} \right) \right)
_{i,j=k+1, \ldots , n}.$$
For $1 \leq \ell \leq k$ and $k+1 \leq j \leq n$ we have
$$\frac{dx_\ell}{dx_j} = (-1)^{k+\ell+1}
\frac{1}{\Delta} \frac{\partial(f_{\varepsilon^\prime 1},
\ldots , f_{\varepsilon^\prime k})}
{\partial(x_1, \ldots, x_{\ell-1}, x_{\ell+1}, \ldots, x_k, x_j)}.$$
Therefore we get
\begin{eqnarray*}
\frac{d}{d x_j}\left( \frac{m_i}{\Delta} \right) & = &
\frac{\partial}{\partial x_j}\left( \frac{m_i}{\Delta} \right) +
\frac{\partial}{\partial x_1}\left( \frac{m_i}{\Delta} \right)
\frac{dx_1}{dx_j}
+ \ldots +
\frac{\partial}{\partial x_k}\left( \frac{m_i}{\Delta} \right)
\frac{dx_k}{dx_j}\\
& = & \frac{1}{\Delta^3} \left( \Delta \frac{\partial m_i}{\partial x_j}
\Delta - m_i
\frac{\partial \Delta}{\partial x_j}\Delta  +(-1)^k \Delta  \frac{\partial
m_i}{\partial x_1}
\frac{\partial(f_{\varepsilon^\prime 1}, \ldots , f_{\varepsilon^\prime
k})}{\partial(x_2, \ldots, x_k, x_j)} \right. \\
& & {} -(-1)^km_i \frac{\partial \Delta}{\partial x_1}
\frac{\partial(f_{\varepsilon^\prime 1}, \ldots , f_{\varepsilon^\prime
k})}{\partial(x_2, \ldots, x_k, x_j)} + \ldots \\
& & {} - \left.  \Delta  \frac{\partial m_i}{\partial x_k}
\frac{\partial(f_{\varepsilon^\prime 1}, \ldots , f_{\varepsilon^\prime
k})}{\partial(x_1, \ldots, x_{k-1}, x_j)} + m_i
\frac{\partial \Delta}{\partial x_k}
\frac{\partial(f_{\varepsilon^\prime 1}, \ldots , f_{\varepsilon^\prime
k})}{\partial(x_1, \ldots, x_{k-1}, x_j)}
\right)\\
& = & \frac{1}{\Delta^3} \left|
\begin{array}{ccccc}
\Delta & \frac{\partial \Delta}{\partial x_1} & \cdots &
\frac{\partial \Delta}{\partial x_k} & \frac{\partial \Delta}{\partial x_j}\\
0 & \frac{\partial f_{\varepsilon^\prime 1}}{\partial x_1} & \cdots &
\frac{\partial f_{\varepsilon^\prime 1}}{\partial x_k} &
\frac{\partial f_{\varepsilon^\prime 1}}{\partial x_j}\\
\vdots & \vdots & \ddots & \vdots & \vdots\\
0 & \frac{\partial f_{\varepsilon^\prime k}}{\partial x_1} & \cdots &
\frac{\partial f_{\varepsilon^\prime k}}{\partial x_k} &
\frac{\partial f_{\varepsilon^\prime k}}{\partial x_j}\\
m_i & \frac{\partial m_i}{\partial x_1}  & \cdots & \frac{\partial
m_i}{\partial x_k} &
\frac{\partial m_i}{\partial x_j}
\end{array} \right|.
\end{eqnarray*}
Proposition~\ref{propHessian} now follows from the following two lemmas.

\begin{lemma}\label{lemmaDet1}
Let $A$ be an invertible $l \times l$-matrix, $B$ an $l \times m$-matrix,
$C$ an $m \times l$-matrix, and $D$ an $m \times m$-matrix. Then
$$\det \left( \begin{array}{cc} A & B \\ C & D \end{array} \right) = \det
A\cdot
\det(D-CA^{-1}B).$$
\end{lemma}

\begin{proof}
Let $E_m$ denote the $m \times m$ identity matrix. Then we have
\begin{eqnarray*}
\det \left( \begin{array}{cc} A & B \\ C & D \end{array} \right) & = & \det
A\cdot \det
\left( \begin{array}{cc} A^{-1} & 0 \\ 0 & E_m \end{array} \right)\cdot
\det \left(
\begin{array}{cc} A & B \\ C & D \end{array} \right) \\
& = & \det A\cdot \det \left( \begin{array}{cc} E_l & A^{-1}B \\ C & D
\end{array} \right) \\
& = & \det A\cdot \det \left( \begin{array}{cc} E_l & A^{-1}B \\ C-CE_l & D -
CA^{-1}B \end{array} \right) \\
& = & \det A\cdot \det \left( \begin{array}{cc} E_l & A^{-1}B \\ 0 & D -
CA^{-1}B \end{array}
\right)\\
& = & \det A\cdot \det (D- CA^{-1}B).
\end{eqnarray*}
\end{proof}

\begin{lemma}\label{lemmaDet2}
Let $A$ be an invertible $l \times l$-matrix, $B$ an $l \times m$-matrix,
$C$ an $m \times l$-matrix, and $D=(d_{ij})$ an $m \times m$-matrix. Let
$c^i$ denote the
$i$-th row of $C$ and $b_j$ the $j$-th column of $B$. Let $H$ be the matrix
$$\left( \left| \begin{array}{cc} A & b_j \\ c^i & d_{ij} \end{array} \right|
\right)
_{i,j=1, \ldots , m}.$$
Then
$$ \det H = (\det A)^{m-1} \left| \begin{array}{cc} A & B \\ C & D
\end{array} \right|.$$
\end{lemma}

\begin{proof}
By Lemma~\ref{lemmaDet1}
$$H = ( \det A (d_{ij}- c^i A^{-1} b_j) )
_{i,j=1, \ldots , m}.$$
Another application of Lemma~\ref{lemmaDet1} yields
$$
\det H  =  (\det A)^m \det (D-CA^{-1}B) = (\det A) ^{m-1} \left|
\begin{array}{cc} A &
B \\ C & D \end{array} \right|.$$
\end{proof}

For any $\varkappa\in\CC^{k+M+n}$ and for any singular point $P$ of the 1-form
$\Omega$ on the (possibly singular) level variety
$V_{\CC,\varkappa}$ (i.e, for a singular point of
$V_{\CC,\varkappa}$ or for a zero of the 1-form on its smooth part),
let $I_P$ be the ideal of $\calO_{P,x}$ generated by the functions
$f_{\varepsilon^\prime i}-\varepsilon_i$ ($i=1,\ldots,k$) and by the
$(k+1)\times(k+1)$-minors of the matrix
$$
\left( \begin{array}{ccc}
\frac{\partial f_{\varepsilon^\prime 1}}{\partial x_1} & \cdots &
\frac{\partial f_{\varepsilon^\prime 1}}{\partial x_n} \\
\vdots & \ddots & \vdots \\
\frac{\partial f_{\varepsilon^\prime k}}{\partial x_1} & \cdots &
\frac{\partial f_{\varepsilon^\prime k}}{\partial x_n}\\
A_1-\alpha_1 & \cdots & A_n-\alpha_n
\end{array} \right).
$$
Let $J_P=\calO_{P,x}/I_P$. From the proof of Theorem~\ref{theoremDim}
it follows that: \newline
1) $\displaystyle
\sum_{P\in V_{\CC,\varkappa}}{\mbox{dim}}_\CC J_P=const=
{\mbox{dim}}_\CC J_0$; \newline
2) if $\varphi_1$, \ldots, $\varphi_L$
are real elements of the
ring $\calO_{\CC^n,0}$ which are representatives of a basis of the
factor-algebra
$J_0$ (considered as analytic functions defined in a common neighbourhood
of the origin in $\CC^n$), then, for $\varkappa$ small
enough, the multigerms of the functions $\varphi_1$, \ldots, $\varphi_L$
are representatives of a basis of the algebra
$J_\varkappa=\oplus_{P\in V_{\CC,\varkappa}}J_P$. \newline
Thus the algebras $J_\varkappa$ form a trivial vector bundle over a
neighbourhood of the origin in $\CC^{k+M+n}$ and the choice of
(real) representatives of the elements of the basis of the algebra $J_0$
fixes a trivialization of this bundle compatible with the real
structure.

As it is usual in similar situations, a quadratic form $Q_\varkappa$
on the algebra $J_\varkappa$ is defined by the formula
$$
Q_\varkappa(\psi_1, \psi_2)=\ell_\varkappa(\psi_1\psi_2),
$$
where $\ell_\varkappa$ is a linear function on the vector space
$J_\varkappa$, $\psi_1\psi_2$
means the product in the algebra $J_\varkappa$.

Let $\varkappa\in(\CC^{k+M+n}, 0)$ be such that $\Delta(P)\ne 0$ for
all singular points $P$ of the 1-form $\Omega_{\vert V_{\CC, \varkappa}}$.
In particular this means that $\varkappa$ is not contained in the discriminant
$\Sigma$ of the deformation and the Hessian $h$ is defined in a neighbourhood
of each singular point $P$. For a singular point $P$ of the 1-form
$\Omega_{\vert V_{\CC, \varkappa}}$, let $\widetilde h_P:= h\cdot\Delta(P)^2$.
For $\psi\in J_\varkappa$, define $\ell_\varkappa(\psi)$ as
$$
\sum_{P\in V_{\CC,\varkappa}}\widetilde\ell_P(\psi),
$$
where $\widetilde\ell_P(\psi)=0$ for $\psi$
from the summand $J_{P^\prime}$
for $P^\prime\ne P$,
$$
\widetilde\ell_P(\psi)=\lim_{\beta\to 0}\sum
\frac{\psi(a_i)}{\widetilde h_P(a_i)}
$$
for $\psi\in J_{P}$, $\beta=\beta_{k+1}dx_{k+1}+\ldots,+\beta_n dx_n$,
the sum is over all simple (i.e., nondegenerate) zeros of the 1-form
$\Omega-\beta$ which emerge from the (generally speaking,
degenerate) zero $P$ for generic $\beta$.
The expression under the limit sign is defined for a generic $\beta$
(i.e., for a generic $\beta$ the form $\Omega-\beta$ has only
nondegenerate zeros). The facts that: it has a finite limit as $\beta$
tends to zero (as an element of $\CC^{n-k}$),
this limit depends analytically on $\varkappa$,
the corresponding quadratic form is nondegenerate, it becomes real
for $\varkappa$ real, and in the last case the signature of the corresponding
real quadratic form is equal to
$$
\sum\limits_{P\in V_{\varepsilon}}{\rm ind}_P\, \omega = {\rm ind}_0\, \omega+
(\chi(V_{\varepsilon})-1)
$$
-- follow from the know properties of the similar objects for the smooth case
(see, e.g., \cite{AGV}, \S~5).

At the moment the required
linear functions $\ell_\varkappa$ (and thus the quadratic forms $Q_\varkappa$)
are defined for $\varkappa$ outside of the discriminant $\Sigma$ and
of the set $\Xi$ of those $\varkappa$ for which on $V_{\CC,\varkappa}$
there exists a zero $P$
of the form $\Omega_{\vert V_{\CC,\varkappa}}$ with $\Delta(P)=0$
(both $\Sigma$ and $\Xi$ are hypersurfaces in $\CC^{k+M+n}$).
We want to show that in
fact this family of linear functions has an analytic extension to these
two subsets as well. Moreover we should control the extension to the
last subset so that the quadratic form $Q_\varkappa$ does not degenerate
there.

To show that the constructed family of linear functions $\ell_\varkappa$
has an analytic extension to the set $\Xi\setminus\Sigma$ of those
$\varkappa\in\CC^{k+M+n}\setminus\Sigma$ for which there exists a zero $P$
of the form $\Omega_{\vert V_{\CC,\varkappa}}$ with $\Delta(P)=0$,
it is sufficient to prove that, for a {\bf nondegenerate} zero $P$
of the form $\Omega_{\vert V_{\CC,\varkappa}}$ with $\Delta(P)=0$,
the linear function $\widetilde\ell_{P^\prime}(\psi)$ has a finite limit
different from zero as $P^\prime$ tends to $P$, where $P^\prime$ is a
singular point of the 1-form $\Omega_{\vert V_{\CC,F(P^\prime)}}$,
$F(P^\prime)\notin\Xi$. In this case $\dim_\CC J_P=\dim_\CC J_{P^\prime}=1$,
the vector spaces $J_P$ and $J_{P^\prime}$ are generated by one element
$\varphi\equiv 1$, $\widetilde\ell_{P^\prime}(1)=
\frac{1}{\widetilde h_{P^\prime}(P^\prime)}$, and it is sufficient to
show that $\widetilde h_{P^\prime}(P^\prime)$ has a finite nonzero
limit as $P^\prime\to P$.  Let $\sigma=\{j_1, \ldots, j_n\}$
be a permutation of the indices $1$, \dots, $n$
such that the Jacobian $\Delta^\prime$ of the functions
$f_1$, \dots, $f_k$ with respect
to the variables $x_{j_1}$, \dots, $x_{j_k}$ at the point $P$
(and thus at all points $P^\prime$ close to $P$) is different
from zero. In this case $x_{j_{k+1}}$, \dots, $x_{j_{n}}$ are
local coordinates on the manifold $V_{\CC,\varkappa}$ at the
point $P$ and thus at all points $P^\prime$ close to $P$ (on the corresponding
level manifold). At a point $P^\prime$ close to $P$ and such that
$\Delta(P^\prime)\ne 0$ as well (and thus where $x_{k+1}$, \dots, $x_n$
are also local coordinates) the Jacobian of the coordinate change
$x_{k+1}$, \dots, $x_n$ $\mapsto$ $x_{j_{k+1}}$, \dots, $x_{j_{n}}$ is equal to
$\mbox{sgn\,}(\sigma)\cdot\frac{\Delta^\prime(P^\prime)}
{\Delta(P^\prime)}$.
Thus the value of the Hessian of the restriction of the form $\Omega$
to the corresponding level manifold in the coordinates $x_{j_{k+1}}$,
\dots, $x_{j_{n}}$ at the point $P^\prime$ is equal to
$$
\frac{h(P^\prime)\Delta(P^\prime)^2}{\Delta^\prime(P^\prime)^2}
$$
and therefore differs from $\widetilde h(P^\prime)$ by a nonzero
analytic factor. This finishes the proof in this case.

To show that the constructed family of linear functions $\ell_\varkappa$
has an analytic extension to the discriminant $\Sigma$ it is sufficient
to prove the following. Let $P$ be a point at which the corresponding level
variety $V_{\CC, \varkappa}$ has a singularity of type $A_1$, the 1-form
$\Omega$ (as a 1-form on $\CC^{n+M+n}$) does not tend to zero, and the
zero hyperplane of it is in general position with respect to the tangent
cone of the variety $V_{\CC, \varkappa}$ at the point $P$. In this case
$\dim_\CC J_P=2$. Let $\widetilde\varkappa$ from a neighbourhood of $\varkappa$
be of the form $(\widetilde\varepsilon_1, \ldots, \widetilde\varepsilon_k,
\varepsilon_{k+1}, \ldots, \varepsilon_{k+M}, \alpha_1, \ldots, \alpha_n)$
(i.e., $\widetilde\varkappa$ differs from $\varkappa$ only by the first $k$
coordinates: the values of the functions $f_1$, \dots, $f_k$) and such that
$\widetilde\varkappa\notin\Sigma$. The 1-form
$\Omega_{\vert V_{\CC,\varkappa}}$ has two nondegenerate zeros
$P_1=P_1(\widetilde\varkappa)$ and $P_2=P_2(\widetilde\varkappa)$.
It is sufficient to show that the linear function
$\widetilde\ell_{P_1}+\widetilde\ell_{P_2}$ has a finite limit as
$\widetilde\varkappa\to\varkappa$.

Without loss of generality one can suppose that $k=1$, $n\ge2$,
$P$ is the origin in $\CC^n$, $f_1=x_1^2+\ldots+x_n^2$,
$\varkappa=0$, $\omega(0)=dx_1$.
The last equation means that $\omega=(1+C_1)dx_1+\ldots+C_ndx_n$,
where $C_i\in{\frak m}$, i.e., $C_i(0)=0$. Let $\widetilde\varkappa=
\varepsilon^2$. As a basis of the algebra $J_P$ (as a vector space) and
thus of the algebra $J_{P_1}\oplus J_{P_2}$ one can take $\varphi_1\equiv 1$
and $\varphi_2=x_1$. For the coordinates of the points $P_1$ and $P_2$
one has $x_1=\pm\varepsilon+ o(\varepsilon)$, $x_i=o(\varepsilon)$ for
$i\ge 2$. Here and further on all series are power series in $\varepsilon$
and thus, for example, $o(\varepsilon)$ means
$a_2\varepsilon^2+~terms~of~higher~degree$. From Proposition~\ref{propHessian}
it follows that $h(P_i)=(-1)^{n-1}(\pm\varepsilon)^{1-n}
+~terms~of~higher~degree$ and thus $\widetilde h(P_i)=
(-1)^{n-1}(\pm\varepsilon)^{3-n}+\ldots$. Therefore
$\widetilde\ell_{P_1}(1)+\widetilde\ell_{P_2}(1)=(-1)^{n-1}
(\varepsilon^{n-3}+(-\varepsilon)^{n-3})+\ldots$,
$\widetilde\ell_{P_1}(x_1)+\widetilde\ell_{P_2}(x_1)=(-1)^{n-1}
(\varepsilon^{n-2}+(-\varepsilon)^{n-2})+\ldots$.
Both expressions have finite limits as $\varepsilon$ tends to zero
(for $n=2$, the terms $\varepsilon^{-1}$ and $(-\varepsilon)^{-1}$
sum up to zero).
\end{proof}

\begin{remarks}
{\bf 1)} To prove Theorem~\ref{theo3} actually it was not necessary
to make precise computations for the $A_1$ case. It is clear that
the components of the linear function $\ell_\varkappa$ (its values
on the basis elements) can have only power asymptotics when $\varepsilon
\to 0$. Thus multiplying the constructed functions by a suitably
high power of the equation of the discriminant we obtain the
required family.
\newline
{\bf 2)} The calculations above for the $A_1$-singularity show that
for $n-k=1$, i.e., for curves, the constructed family of quadratic
forms does not degenerate anywhere (including the discriminant) and,
in particular, the quadratic form $Q_0$ (defined on the algebra $J_0$) is
nondegenerate and its signature is equal to the same expression (\ref{eq2}).
\newline
{\bf 3)} In \cite{Szafraniec01} there was proved a somewhat similar
statement which in our terms can be considered as a particular case for
the 1-form $\omega$ being the differential $df_{k+1}$ of a function.
However the family of quadratic forms defined there on the target space
of the map $f:\CC^n\to\CC^k$ could degenerate at points $\varepsilon\in
\CC^k$ for which there exists a zero $P$ of the 1-form $\omega$ on
the level manifold $V_{\CC,\varepsilon}$ with $\Delta(P)=0$ (and maybe
somewhere else). (By the way, in some particular cases all points
$\varepsilon$ from the target $\CC^k$ could possess this property.)
Thus one can say that our result gives a certain improvement of the
one from \cite{Szafraniec01} for $\omega=df_{k+1}$ as well.
\end{remarks}

\bigskip
\noindent Institut f\"{u}r Mathematik, Universit\"{a}t Hannover, \\
Postfach 6009, D-30060 Hannover, Germany \\
E-mail: ebeling@math.uni-hannover.de\\

\medskip
\noindent Department of Mathematics and Mechanics,\\
Moscow State University\\
Moscow, 119899, Russia\\
E-mail: sabir@mccme.ru

\end{document}